\documentclass[11pt, reqno]{amsart}
\usepackage{graphicx}
\usepackage{amssymb}

\newtheorem{theorem}{Theorem}[section]
\newtheorem{lemma}[theorem]{Lemma}
\newtheorem{corollary}[theorem]{Corollary}
\newtheorem{proposition}[theorem]{Proposition}

\theoremstyle{definition}

\def\ord{\operatorname{ord}}

\def\beq {\begin{equation}}
\def\endq {\end{equation}}
\def\bs{\boldsymbol}
\def\e{\varepsilon}
\def\d{\delta}

\def\beq {\begin{equation}}
\def\endeq {\end{equation}}
\newcommand{\ignore}[1]{#1}

\renewcommand{\epsilon}{\varepsilon}

\begin{document}
\title{The Uniform Primality Conjecture for the Twisted Fermat Cubic\ }
\thanks{The research of the second author was
supported by a grant from the Thai government.}

\author{Graham Everest, Ouamporn Phuksuwan and Shaun Stevens}
\address{School of Mathematics, University of East Anglia,
Norwich NR4 7TJ, UK}
\dedicatory{To Jupiter, bringer of joy.}
\email{g.everest@uea.ac.uk}
\email{o.phuksuwan@uea.ac.uk}
\email{s.stevens@uea.ac.uk}

\begin{abstract}
On the twisted Fermat cubic, an elliptic divisibility
sequence arises as the sequence of denominators of the
multiples of a single rational point. We prove that the number of prime
terms in the sequence is uniformly bounded. When the rational point
is the image of another rational point under a certain $3$-isogeny,
all terms beyond the
first fail to be primes.

\end{abstract}

\maketitle


\section{Introduction}

There abides a deep fascination with the subject of
prime terms in naturally occurring integer sequences. The fanfare
which accompanies the discovery of each new Mersenne prime \cite{caldwell1}
bears testimony,
as
does the search for factorizations of the Fermat and Fibonacci numbers
\cite{caldwell2,keller}.
Elliptic curves \cite{cassels,aec,Silverman-advance} give rise naturally to rapidly growing sequences of integers,
called {\it elliptic divisibility sequences},
generated from some nontorsion point via the geometric addition law.
The first systematic study of elliptic divisibility sequences was given
majestically by Morgan Ward in \cite{ward}. Recent times have witnessed a
great deal of activity (the following represents a sample) because of their intrinsic interest
\cite{emw,pi1,pi3,pi2,sileds1,sileds3,silstang,sileds2,streng} as well as in applications to
Logic \cite{cz,ee,bjorn,ss,ssbook}
and
Cryptography \cite{ls,shipsey-thesis}.
The Chudnovsky brothers~\cite{MR88h:11094} first proposed studying primality in
elliptic divisibility sequences arising from
elliptic curves in Weierstrass form. This was
continued in \cite{eds,pe,primeds} with the following outcome: a heuristic
argument,
together with lots of computational evidence and proofs in some special cases,
suggests that any elliptic divisibility sequence contains only
finitely many prime terms. This has become known as the {\it primality
conjecture} for elliptic curves in Weierstrass form. In \cite[Theorem 4.1]{primeds}
the same conjecture was proved for curves in homogeneous form. It is with such curves that
our business lies in this paper.

Let $C$ denote an elliptic curve of the form
\begin{equation}\label{defC}
C : U^3 + V^3 = m,
\end{equation}
where $m$ is a nonzero integer. It is sometimes said that~$C$ is a twist of
the Fermat cubic. The set of rational (or real) points form a group under the
chord and tangent method: the (projective) point at infinity with homogeneous coordinates $(1:-1:0)$ is the
identity, and inversion is given by reflection in the line $U=V$.
Given a nontorsion point $R \in C(\mathbb{Q})$, write, in lowest terms,
\begin{center}
$nR = \left( \dfrac{U_n}{W_n}, \dfrac{V_n}{W_n} \right).$
\end{center}
The sequence $(W_n)$ is a (strong) divisibility sequence, see
Proposition~\ref{strongdivseq}.
Our starting point is the result
stated earlier, that only finitely many terms
$W_n$ are prime \cite[Theorem 4.1]{primeds}.
A related result \cite[Theorem 1.2]{spag} states that $W_n$ has a primitive divisor
for all~$n>1$: results from \cite{spag} will be employed in the proofs of
our theorems.

To focus now upon the question we study: a stronger form of the primality conjecture has been
proposed, assuming that the Weierstrass
equation is
in minimal form, predicting
that the number of prime terms in an elliptic divisibility sequence is
uniformly bounded (independently of both the curve and the point). Unsurprisingly, this conjecture is known as
the {\it uniform primality conjecture} for curves in Weierstrass form.
A conditional proof is known, under Lang's Conjecture,
assuming the generating point is the image of a rational point under
a nontrivial isogeny \cite[Theorem 1.4]{eims}. This raises the natural question as to whether
a uniform primality result holds for curves in homogeneous form and it is this
question we address. The first of our two results is an unconditional result to that
effect, and will be stated now.
The proof will appear to be quite short although, in fact, it relies upon two substantial
results from \cite{spag} and \cite{MR2468477}.

\begin{theorem}\label{upc}Assuming $m\in \mathbb Z$ is cube-free, the
number of prime (power) terms~$W_n$ is uniformly bounded. Moreover, assuming the $ABC$ conjecture holds in~$\mathbb Z$, there
is a uniform upper bound on the index~$n$ such that~$W_n$ is
a prime power.
\end{theorem}

The assumption about $m$ being cube-free is necessary given the homogeneous nature
of the curve~$C$, since one could scale any given elliptic divisibility sequence to obtain as many prime terms as desired.
This will be explained in detail in subsection~\ref{cubefree}.
It is also probable that a uniform bound on the index~$n$ would follow from the Hall-Lang conjecture,
as in recent work of Mah\'e~\cite{valery}.

Our second theorem relies upon the existence of a rational isogeny in line
with work for curves in Weierstrass form \cite{pe,primeds}.
There is a $3$-isogeny, detailed in \cite{rande},
$\tau:C'\rightarrow C$ from the elliptic curve
\begin{equation}\label{defC'}C': U'V'(U'+V')=m.
\end{equation}

\begin{theorem}\label{expupc}
Suppose $m\in \mathbb Z$ is cube-free and $R\in C(\mathbb Q)$ is the image of a point in $C'(\mathbb Q)$ under $\tau$.
Then $W_n$
is not a prime power, for all $n > 1$. Further, the number of terms~$W_n$ with at most two
distinct prime divisors is uniformly bounded.
\end{theorem}

The conclusion is clearly best possible and what is remarkable is the small size of the bound~$1$; this is certainly much stronger than any corresponding result for elliptic curves in Weierstrass form.

The proofs of the two theorems rely upon the same principle, which will
be articulated in section~\ref{proveit} below. This is preceded by
a section discussing possible improvements to our results and is followed by proofs of the
theorems in order.

We thank Tony Flatters and Val\'ery Mah\'e for their excellent comments and suggestions on a preliminary
version of this paper. In particular, the latter pointed out that a (conditional) bound on the index could be
obtained in Theorem~\ref{upc}.

\section{Better-best?}

\subsection{With isogeny}
Taking our two theorems in reverse over, Theorem \ref{expupc} cannot be improved.
But to how large a class of examples
does it apply? If~$R$ is $3$ times
another rational point then it will be such an image~\cite{rande}. A short analysis
is included to shed light in greater generality.
In truth, throughout the paper, we will not work
explicitly with the isogeny mentioned in the statement of Theorem~\ref{expupc}.
Instead, our work relies upon a bi-rational map between~$C$ and the Mordell curve
\begin{equation}\label{defE}
E: Y^2 = X^3 - 432m^2,
\end{equation}
where $P=(X,Y)$ corresponds to $R=(U,V)$ under the bi-rational transformation given by
\begin{equation}\label{transform}
\begin{array}{lr}
X = \dfrac{2^{2}3m}{U+V},& Y = \dfrac{2^{2}3^{2}m(U-V)}{U+V},\\[15pt]
U = \dfrac{2^{2}3^{2}m + Y}{6X},& V = \dfrac{2^{2}3^{2}m - Y}{6X}.
\end{array}
\end{equation}

There is a 3-isogeny $\sigma:E'\rightarrow E$ onto the curve $E$ from the elliptic curve
$$
E' : y^2 = x^3 + 16m^2,
$$
given by
\begin{equation}\label{isog}
\sigma (x, y) = (X, Y) = \left( x + \dfrac{64m^2}{x^2}, \dfrac{y(y+12m)(y-12m)}{(y+4m)(y-4m)}\right).
\end{equation}
It is with this isogeny that we work. The isogeny condition stated in Theorem~\ref{expupc} is
equivalent to the point $P\in E(\mathbb Q)$ being equal to $\sigma(P')$ for some point $P' \in E'(\mathbb Q)$.
Table~\ref{m<100} shows the $22$ values $m < 100$ when this occurs for a point~$P$ which is a generator for the group~$E(\mathbb Q)$.
(Thus all points in~$E(\mathbb Q)$ are in the image of the isogeny.)
In total, there are
$42$ rank one curves with~$m<100$. We
conclude that, at least for small values of~$m$, the hypothesis of
Theorem~\ref{expupc} is not infrequently met.

\begin{table}[!ht]
\caption{Rank 1 curves satisfying Theorem~\ref{expupc} with $m<100$}
\label{m<100}
  \centering
  \begin{tabular}{|c|c|c|}
   \hline
     $m$&$P$&$P'$\\
     \hline
6& [28, 80]& [-8, 8] \\
12&[52, 280]&[-12, 24]\\
15& [49, 143] & [-15, 15]\\
20& [84, 648]& [-16, 48]\\
22& [553/9, 4085/27]& [33, 209] \\
33& [97, 665]& [-24, 60] \\
34& [273, 4455]& [-16, 120] \\
42 & [172, 2080] & [-24, 120]\\
50 & [2716/9, 138736/27] & [24, 232]\\
51& [5473/36, 333935/216] & [144, 1740]\\
58 & [3148/9, 173600/27] & [-24, 200]\\
68& [55156/225, 12043304/3375 ]& [240, 3728] \\
69& [553, 12925]& [-23, 253]\\
70& [156, 1296] & [-40, 120]\\
75& [601, 14651]& [-24, 276]\\
78 & [217, 2755] & [-39, 195]\\
84 &[148, 440] & [-48, 48]\\
85& [18361/9, 2487509/27] & [-15, 335] \\
87&[1029841/1225, 1042214111/42875]& [840, 24348]\\
90& [364, 6688] & [-36, 288]\\
92 & [2548/9, 117800/27] & [48, 496]\\
94 &[62511752209/2480625,\hfill\ &\\
&\ \hfill15629405421521177/3906984375]& [25200, 4000376]\\
\hline

  \end{tabular}
\end{table}

\subsection{Without isogeny}
Turning now to Theorem \ref{upc}, we will draw upon a result of Ingram \cite{MR2468477}
bounding the number of integral multiples of a given integral point. In its
most general form, this involves a non-explicit constant. In his paper, he
gives a very strong (actually best possible) explicit constant in the case
of congruent number curves. It seems a good bet that his analysis could be
re-worked in our context to obtain an explicit (perhaps even best possible)
form of Theorem~\ref{upc}.
This analysis was done in certain cases by the second author in her thesis~\cite{opthesis}:
a special case of \cite[Theorem~1.0.4]{opthesis} is recorded here.

\begin{theorem}Assume $m \equiv \pm 1, \pm 3, \pm 4$ mod $9$ is cube-free.
Assuming $P\in E(\mathbb Q)$ has $\gcd (x(P),m)=1$ and~$2P$ and~$3P$ are not integral
then at most one value of $n>1$ yields a prime power term~$W_n$.
\end{theorem}

Of course the most desirable goal is to be to be able
to compute all of the prime power terms in a given sequence, as permitted by
Theorem~\ref{expupc}. It is not clear whether Theorem~\ref{upc} could be
improved to allow that as an outcome. However, as part of the proof of Theorem~\ref{expupc},
we do prove the following unconditional result:

\begin{lemma}\label{lem:2<n<23}
Suppose $m\in \mathbb Z$ is cube-free. Then $W_n$
is not a prime power, for $2 < n \le 22$.
\end{lemma}

\subsection{Daylight}
There will always be daylight between out theorems as our final remark
demonstrates. It seems likely that without the isogeny hypothesis, there
are infinitely many values of~$m$ and generating points~$R$ yielding a prime term~$W_2$.
For, suppose $R = [u,v] \in C(\mathbb{Q})$
is an integral point of infinite order.
The point $2R$ can be expressed in the form
$$
2R = \left(\frac{-2vu^3 - v^4}{u^3 - v^3}, \frac{u^4 + 2v^3u}{u^3 - v^3}\right).
$$
Suppose $u - v = 1$. Then
\begin{center}
$u^3 - v^3 = (u-v)(u^2 + uv + v^2)  = 3u^2 - 3u + 1$.
\end{center}
Applying the Bateman-Horn conjecture \cite{MR0148632} to the polynomial
\begin{center}
$f(u) := 3u^2 - 3u + 1$
\end{center}
implies that $f(u)$ is prime for infinitely many positive integers $u$.
By a result of Erd\~os~\cite{erdos}, the corresponding value of $m=u^3+(u-1)^3$ will be cube-free for infinitely
many of these values. To show
that the prime so constructed does not cancel with the numerator, compute the resultants of the two polynomials
$u^2+uv+v^2$ and $2u^3+v^3$ in each variable: they are $9u^6$ and $9v^6$ so the only non-trivial
common division is by a divisor of~$9$. But $3u^2 - 3u + 1$ is coprime to $3$ so there can be no cancellation.
This style of argument will be used
repeatedly in the proof of Theorem~\ref{expupc}.

\section{Proving the Primality Conjecture}\label{proveit}

This section is elementary although quite technical because of the case-by-case
nature of the proofs. It is the statements of the results which
provide important input for sections \ref{proofofupc} and~\ref{proof1.2}, rather than their proofs.
Readers whose taste is not
excited by the
details could take note of Lemma~\ref{cancel}, Proposition~\ref{strongdivseq} and
subsection~\ref{cubefree} then skip directly to sections~\ref{proofofupc} and~\ref{proof1.2}.

\subsection{Constraints to cancellation} From the bi-rational transformation \eqref{transform}, we have
\begin{equation}\label{UnWn}
\frac{U_n}{W_n} = \frac{2^{2}3^{2}mB_{n}^{3} + C_n}{6A_{n}B_n},
\end{equation}
and
\begin{equation}\label{VnWn}
\frac{V_n}{W_n} = \frac{2^{2}3^{2}mB_{n}^{3} - C_n}{6A_{n}B_n},
\end{equation}
where $$ nR=\left(\dfrac{U_n}{W_n},\dfrac{V_n}{W_n}\right) \mbox{ and } nP =
\left( \dfrac{A_n}{B_{n}^{2}}, \dfrac{C_n}{B_{n}^{3}} \right)$$
are written in lowest terms.

\begin{lemma}\label{cancel}Any cancellation in the fractions \eqref{UnWn} and \eqref{VnWn}
comes from the term~$6A_n$ and divides~$72m$.
\end{lemma}

Since $A_n$ and $B_n$ are coprime, the following important principle arises as a consequence of Lemma \ref{cancel}:

\vskip0.1in

\begin{center}
\begin{minipage}[c]{4.6in}
{\bf
We can be certain that $W_n$ possesses at least two coprime factors if two
conditions are met. The first is that~$B_n>1$ and the second is that the term~$6A_n$ is
not completely cancelled.}
\end{minipage}
\end{center}

\vskip0.1in

\begin{proof}[Proof of Lemma \ref{cancel}]
Consider the fraction on the right-hand side of \eqref{UnWn}.
Let $d = p^r$ be a common factor of
$(2^{2}3^{2}mB_{n}^{3} + C_n)$ and $6A_nB_n$ with $p$ a prime number and $r \in \mathbb{N}$ the highest order of $p$
dividing both terms. If $d' := \gcd(d, B_n) \neq 1$, then $d' \mid  d \mid  (2^{2}3^{2}mB_{n}^{3} + C_n)$ implies $d' \mid  C_n$ ,
which contradicts the fact that $B_n$ and $C_n$ are coprime. Thus $\gcd(d, B_n) = 1$, so that $d$ comes from the term $6A_n$.

Notice, moreover, that any cancellation in the right-hand side of \eqref{UnWn} and \eqref{VnWn} is the same
in each term.
This is because both terms have the same denominators as their cubes sum to an integer.
Hence $d$ has to divide both $(2^{2}3^{2}mB_{n}^{3} + C_n)$  and $(2^{2}3^{2}mB_{n}^{3} - C_n)$, so that $d \mid  72m$.
Thus any cancellation of the fractions \eqref{UnWn} and \eqref{VnWn} divides $72m$.
\end{proof}

We will need to be even be more precise about the cancellation that occurs.
Write $d_n$ for the cancellation in the fraction~\eqref{VnWn}.

\begin{corollary}\label{cor:cancel}
We have
$$
\ord_p(d_n)=\begin{cases} \ord_p(m) &\hbox{ if }p\mid A_n,\ p>3; \\
\ord_p(A_n)+1 &\hbox{ if }p\mid A_n,\ p\le 3; \\
0&\hbox{ if }p\nmid A_n. \end{cases}
$$
\end{corollary}

Note in particular that if $p\le 3$ and $p\mid A_n$ then $p$ is completely cancelled. 
Thus for primes $p\le 3$, we have $p\mid W_n$ if and only if $p\nmid A_n$. For primes 
greater than $3$, the cancellation is exactly as large as allowed by Lemma~\ref{cancel}; 
therefore, whenever $p>3$ divides $A_n$, we have
\begin{equation}\label{alter}
\ord_p(W_n) = \ord_p(A_n)-\ord_p(m).
\end{equation}

\begin{proof}
Suppose first that $p\nmid A_n$. By Lemma~\ref{cancel}, we have $\ord_p(d_n)=0$ except 
possibly if $p\le 3$ so suppose we are in this situation. The defining equation~\eqref{defE} yields
\begin{equation}\label{eqnABC}
C_n^2=A_n^3-432m^2B_n^6,
\end{equation}
whence $p\nmid C_n$ so $p$ cannot cancel in~\eqref{VnWn}.

Now suppose $p\mid A_n$ and $p\ge 3$. From~\eqref{eqnABC} we find $$\ord_p(A_n)\ge \ord_p(C_n)= \ord_p(m),$$ 
since $p\mid A_n$, and it follows that $\ord_p(d_n)=\ord_p(m)$.

Finally, suppose $p\mid A_n$ and $p\le 3$, and put $\d=\ord_p(\gcd(p,m))$. The arguments are 
fiddly but elementary so we only sketch them, beginning with the slightly simpler case~$p=3$.

If $\ord_3(A_n)> 1+\d$ then, 
comparing $3$-adic valuations, we see that
equation~\eqref{eqnABC} is impossible to solve.
Hence $1\le \ord_3(A_n)\le 2$ and we also find that $\ord_3(C_n)\ge \ord_3(A_n)+1$. Thus the numerator of~\eqref{VnWn} has $3$-adic order at least~$\ord_3(A_n)+1$ whereas the $3$-adic order of the denominator is exactly~$\ord_3(A_n)+1$, giving $\ord_3(d_n)=\ord_3(A_n)+1$.

\medskip

Now assume that $p=2$. If $\ord_2(A_n)> 2+\d$ then, dividing~\eqref{eqnABC} through by $2^4 3^2m^2$ yields a congruence $x^2 \equiv -1\pmod4$, which is impossible. Moreover, $\ord_2(A_n)=1$ is impossible by comparing $2$-adic valuations in~\eqref{eqnABC}.

If $\ord_2(A_n)=2+\d$ then $\ord_2(C_n)=2+\ord_2(m)$ so
the numerator of~\eqref{VnWn} is $2^{2+\ord_2(m)}$ times the sum of two odd numbers. Thus the $2$-adic order of the numerator is at least $3+\ord_2(m)$ while that of the denominator is exactly~$3+\d$ so
$\ord_2(d_n)=\ord_2(A_n)+1$.

This leaves the case $\ord_2(A_n)=2$ with $2\mid m$. Here $\ord_2(C_n)>2$ so the $2$-adic order of the numerator of~\eqref{VnWn} is at least $3$, while that of the denominator is exactly~$3$, again giving $\ord_2(d_n)=\ord_2(A_n)+1$.
\end{proof}

\subsection{$\bs W$ is a divisibility sequence}
The aim of this subsection is a proof of the following:
\begin{proposition}\label{strongdivseq}
The sequence $(W_n)$ is a strong divisibility sequence: in other words, for all $r, n \in \mathbb N$,
$$
\gcd(W_r,W_n)=W_{\gcd(r,n)}.
$$
In particular, $W_n\mid W_{nk}$, for all $n,k \in \mathbb N$.
\end{proposition}

The divisibility property
will be used repeatedly in this paper but
we cannot find it proved explicitly in the literature.
The strong divisibility property
will be used only in subsection~\ref{cubefree}. However, it is
a natural property and is proved with little more effort.

We admit that the proof we give lacks elegance. This is due to the evil influence of the primes~$2$
and~$3$ and, to a lesser extent, those dividing~$m$. However, the proof is completely self-contained
and uses only elementary methods.

A more sophisticated proof of Proposition~\ref{strongdivseq} uses formal groups~\cite[Chapter IV]{aec}.
Writing $w=U/V$ and $z=1/V$ in (\ref{defC}) yields
the equation
$$w^3=mz^3-1.
$$
Now the binomial theorem yields a power series
$$w=\sum_{n=0}^{\infty}a_nz^{3n} \mbox{ with } a_n \in \mathbb Q,
$$
which converges $p$-adically for all primes $p$ with $p\mid z$ (only the case $p=3$
is at all tricky). Applying the geometric
group law on the points $(z_1,w_1)$ and $(z_2,w_2)$ now yields a power series
$F(z_1,z_2)$ -- the {\it formal group} of our elliptic curve. The proof of
Proposition~\ref{strongdivseq} 
follows from the standard properties of formal groups~\cite[Chapter IV, VII]{aec}:
in particular, the filtration into subgroups~$C_r$,~$r\in \mathbb N$, defined by
$|U|_p\ge p^r$.
However arguing as above requires quite a lot of explanation and checking and 
it would not shorten the paper. 
For example, the statement of Proposition~\ref{strongdivseq} per se does not appear in~\cite{aec}.
Instead, the filtration referred to yields directly the relation~(\ref{W-ordthing}) below. 

\begin{proof}[Proof of Proposition~\ref{strongdivseq}]
The divisibility property follows from
the following relation: for any prime~$p$ and $k\in\mathbb N$, if $p\mid W_n$ then
\begin{equation}\label{W-ordthing}
\ord_p(W_{nk})=\ord_p(W_n)+\ord_p(k).
\end{equation}

Much of the spade work here is supplied by the following two relations. Firstly, from \cite[Lemma 3.2]{spag},
for~$p>3$, if $\ord_p(A_n)>0$ and $3\nmid k$ then
\begin{equation}\label{ordthing}
\ord_p(A_{nk})=\ord_p(A_n)+\ord_p(k).
\end{equation}
Secondly, for all primes $p$, if $\ord_p(B_n)>0$ and $k\in\mathbb Z$ then
\begin{equation}\label{B-ordthing}
\ord_p(B_{nk})=\ord_p(B_n)+\ord_p(k).
\end{equation}
The second relation arises from a local analysis using the formal
group for the Weierstrass equation and follows from~\cite[Chapter VII]{aec}.
A much more general proof, valid over a
Dedekind domain, is given in~\cite[Proposition 1]{pi2}.

Note that~\eqref{B-ordthing} implies that $B=(B_n)$ is a
strong divisibility sequence. Further work is needed to take
us from (\ref{W-ordthing}) to the same conclusion for the sequence~$W=(W_n)$: this will be supplied at the end.

Assume $p$ denotes a prime with~$p\mid W_n$. We will show~\eqref{W-ordthing} for all $k\in \mathbb N$. There are
three cases in total:
\begin{itemize}
\item[(a)] \ $p\mid B_n$,
\item[(b)] \ $p\nmid B_n,\ p>3$,
\item[(c)] \ $p\nmid B_n,\ p\le 3$.
\end{itemize}

\smallskip

{\bf Case (a)} No cancellation occurs in~\eqref{VnWn} so $\ord_p(W_n)=\ord_p(B_n)$ if $p\nmid 6$ or
$1+\ord_p(B_n)$ if $p \mid 6$. Then~\eqref{W-ordthing} follows from~\eqref{B-ordthing}, the corresponding result
for the sequence~$B=(B_n)$.

\smallskip

{\bf Case (b)} Since $p\mid W_n$, the two conditions given imply $p\mid A_n$. 
When $3\nmid k$, we have $p\mid A_{nk}$ by~\eqref{ordthing} so we can apply~\eqref{alter} with $nk$ replacing $n$. Then
\begin{align*}
\ord_p(W_{nk})&=\ord_p(A_{nk}) -\ord_p(m) &&\mbox{by \eqref{alter},}\\
&=\ord_p(A_n)+\ord_p(k)-\ord_p(m) &&\mbox{by \eqref{ordthing},}\\
&=\ord_p(W_n)+\ord_p(k) &&\mbox{by \eqref{alter}}.
\end{align*}
When $3\mid k$, we use the
triplication law~\cite[(8)]{spag}, repeated here for convenience:
\begin{equation}\label{triple}
x(3Q)=
\frac{x^9(Q)+2^93^4x^6(Q)m^2+2^{12}3^7x^3(Q)m^4-2^{18}3^9m^6}{9x^2(Q)(x^3(Q)-2^63^3m^2)^2}.
\end{equation}
It follows from~\eqref{triple} and~\eqref{alter} that
\begin{equation}\label{ordtriple}
\ord_p(B_{3n})=\ord_p(A_n)-\ord_p(m)=\ord_p(W_n).
\end{equation}
Since $3n\mid nk$ we deduce from~\eqref{B-ordthing} that $p\mid B_{nk}$. Therefore
\begin{align*}
\ord_p(W_{nk})&=\ord_p(B_{nk}) &&\mbox{since $p>3$,}\\
&=\ord_p(B_{3n})+\ord_p(k/3) &&\mbox{since $p\mid B_{3n}$,}\\
&=\ord_p(B_{3n})+\ord_p(k) &&\mbox{since $p>3$,}\\
&=\ord_p(W_n)+\ord_p(k) &&\mbox{by~\eqref{ordtriple}}.
\end{align*}

\smallskip

{\bf Case (c)} Since $p\mid W_n$ and $p\le 3$, we have $p\nmid A_n$ and, as there is no cancellation, $\ord_p(W_n)=1$. Also,
since $p\nmid A_n$ it follows that $p\nmid A_{nk}$ for all $k \in \mathbb N$: this is because
the point $nP$ has good reduction at~$p$ so $nkP$ has good reduction mod~$p$
as the set of points with good reduction form a group 
and it follows that $nkP$ cannot reduce mod~$p$ to a
point with zero $x$-coordinate. Thus there is no cancellation in~\eqref{VnWn} when $n$ is replaced by $nk$. Moreover, as $p\nmid B_n$ also, the reduction of $nP$ is a non-identity point and, since $E$ has additive reduction modulo $p$, the point $nkP$ reduces to the identity if and only if $p\mid k$. (Note that we are reducing the curve $E$, not its minimal form.)

There are now two possibilities: If $p\nmid k$ then $p\nmid B_{nk}$, in which case $\ord_p(W_{nk})=1=\ord_p(W_n)+\ord_p(k)$. Otherwise $p\mid k$ and $p\mid B_{nk}$; it is easily checked by an explicit calculation with the doubling and tripling formulae (see~\eqref{triple} for the latter) that $\ord_p(B_{np})=1$ so we get
\begin{align*}
\ord_p(W_{nk})& =\ord_p(B_{nk})+1 \\
&=\ord_p(B_{np}) + \ord_p(k/p)+1 &&\mbox{by~\eqref{B-ordthing}},\\
&=1+\ord_p(k) \\
&=\ord_p(W_n)+\ord_p(k).
\end{align*}

Now we have proved~\eqref{W-ordthing}, the strong divisibility property comes from the following claim:

{\bf Claim} If a prime $p$ divides $\gcd(W_r,W_n)$ then $p\mid W_{\gcd(r,n)}$.

To see that this suffices, let $d=\gcd(r,n)$ and suppose $r=dk,n=dl$. For any
prime~$p$, one of $\ord_p(k)$ and $\ord_p(l)$ is~$0$ because~$k$ and~$l$ are
coprime. From~\eqref{W-ordthing}, the exact power of~$p$ dividing $\gcd(W_{dk},W_{dl})$ is $\ord_p(W_d)$ as desired.

\smallskip

It remains to prove the claim. First let~$p$ be a prime not dividing~$6m$. Note that $p\mid W_n$ if and only if $p\mid A_nB_n$ in this case (since no cancellation can occur), which is equivalent to $nP$ reducing mod $p$ to an element of the subgroup $G$ of $E(\overline{\mathbb F}_p)$ generated by a $3$-torsion point~$T$ with $x(T)=0$. If $rP$ and $nP$ both reduce to an element of $G$ then so does any integer linear combination; in particular, so does $\gcd(r,n)P$ and we deduce that $p\mid W_{\gcd(r,n)}$.

Now suppose $p>3$ is a prime dividing $m$ such that $p\mid\gcd(W_r,W_n)$. There are three possibilities, and in each case we will show that $p\mid W_{\gcd(r,n)}$:
\begin{enumerate}
\item[(i)] $p\mid B_r$ and $p\mid B_n$, in which case $p\mid B_{\gcd(r,n)}$, since $B$ is a strong divisibility sequence. Thus $p\mid W_{\gcd(r,n)}$, since there is no cancellation.

\smallskip
\item[(ii)] $p\mid B_r$ and $p\mid A_n$, in which case $nP$ is a point of bad reduction, and $P$ is also. Thus $p\mid A_1$ and, from the analysis in case (b) above, we see that
$3\mid r$, $3\nmid n$. Then, from~\eqref{ordtriple}, we have $p\mid B_{3n}$ so, from strong divisibility, $p\mid B_{\gcd(r,3n)}=B_{3\gcd(r,n)}$. Moreover, since $3\nmid\gcd(r,n)$, we have $p\mid A_{\gcd(r,n)}$ so, from~\eqref{ordtriple} again,
$$
\ord_p(W_{\gcd(r,n)}) = \ord_p(B_{3\gcd(r,n)}) > 0.
$$

\smallskip
\item[(iii)] $p\mid A_r$ and $p\mid A_n$, in which case $P$ is again a point of bad reduction, $3\nmid r$, $3\nmid n$, so $p\mid A_{\gcd(r,n)}$. From~\eqref{ordtriple}, we have $p\mid\gcd(B_{3r},B_{3n})=B_{3\gcd(r,n)}$ so that $\ord_p(W_{\gcd(r,n)}) = \ord_p(B_{3\gcd(r,n)}) > 0$.
\end{enumerate}

Finally, suppose $p\le 3$. Since $p\mid W_n$ we have $p\nmid A_n$ so $nP$ is a point of good reduction modulo $p$. The same is true of $rP$ and so also of $\gcd(r,n)P$. Thus $p\nmid A_{\gcd(r,n)}$ and we deduce from Corollary~\ref{cor:cancel} that $p\mid W_{\gcd(r,n)}$.

This completes the proof of the claim, and of Proposition~\ref{strongdivseq}.
\end{proof}

\subsection{Cube-free $\bs m$}\label{cubefree} We conclude with a remark justifying the assumption
in both our theorems that~$m$ is cube-free.

The key remark needed follows from Proposition~\ref{strongdivseq}: if $\gcd(m,n)=1$ then $\gcd(W_m,W_n)=1$.
From \cite[Theorem 1.2]{spag}, the terms $W_n$ with~$n>1$ all see a primitive prime
divisor: in other words, a prime~$p\mid W_n$
which has not appeared in an earlier term. Let~$S$ denote a finite set of primes.
For each~$l\in S$, write $W_l=w_lp_l^{e_l}$
with $w_l, e_l \in \mathbb N$ and $p_l\nmid w_l$ equal to a primitive prime divisor. Now write
$$
M=\prod_{l\in S}w_l
$$
and rescale equation~\eqref{defC} by multiplication with~$M^3$. If~$W_l'=W_l/\gcd(W_l,M)$
denotes the resulting denominator then~$W_l'$ is a prime power by the starting remark. Therefore, by expanding~$S$ arbitrarily, we see that,
at the cost only of multiplying~$m$ by a cube, we may produce elliptic divisibility
sequences with arbitrarily large
numbers of prime power terms.

\section{Proof of Theorem \ref{upc}}\label{proofofupc}

\begin{proof}[Proof of Theorem \ref{upc}]
We will prove firstly that, for $n>12$, the term $A_n$ in the denominator is not
completely cancelled. From~\cite[Theorem 2.3]{spag}, all
terms~$A_n$ with~$n>12$ possess a primitive prime divisor~$p_n$.
Then~$p_n\nmid A_1$, because~$p_n$ is a primitive prime divisor, so it follows that $P$ is a point
of good reduction for~$p_n$. Thus all multiples of~$P$, in particular~$nP$, are points of good reduction. However, for primes dividing $6m$, the point on the reduced
curve with $x$-coordinate zero is a point of bad reduction so it follows that $p_n\nmid 6m$.
In particular,~$p_n$ is not cancelled because any cancellation divides~$72m$.

Secondly, from~\cite[Theorem 1]{MR2468477}, with the notation used there,
there is a uniform constant $C$ such that
$n > CM(P)^{16}$ forces $B_n > 1$, except for at most one value of $n$. The quantity $M(P)$ is related to
the Tamagawa number. Since the Mordell curve $E$ has integral $j$-invariant, along the same lines as in
\cite{MR2468477}, $E$ always has $M(P) \leq 6$. It follows
that the number of prime power terms $W_n$ is bounded by~$C6^{16}+1$, a uniform constant.
\ignore{

Finally, assume the $ABC$ conjecture holds for~$\mathbb Z$. As above we always get a prime factor from $A_n$,
for $n>12$, so we need only show that $B_n>1$ for all $n$ greater than some uniform bound. For this we
will use the theory of heights, which will also be essential for the explicit bound in Theorem~\ref{expupc}.

Recall that the \emph{(na\"ive) Weil height} of the point $P$ is
$$
h(P)=h(x(P))=\log\max\left\{|A_1|,B_1^2\right\}.
$$
On the other hand the \emph{canonical height} $\hat h(P)$ is given by
$$
\hat h(P) = \lim_{n\to\infty}\frac{h\left(2^nP\right)}{4^n}.
$$
Silverman gives explicit upper and lower
bounds for the difference between the
Weil height and the canonical height for curves in short
Weierstrass form. Note that our heights are twice those in~\cite{MR1035944} so the inequalities are multiplied by~$2$.

\begin{theorem}[{\cite[Remark 1.2]{MR1035944}}]\label{WeilCan}
Given an elliptic curve in short Weierstrass form,
$$
E/\mathbb{Q}: y^2 = x^3 + ax + b,
$$
and $Q\in E(\mathbb Q)$, we have
\begin{equation}\label{diff}
-\frac{1}{6}h(j) - \frac{1}{6}h(\Delta) - 2.14 \leq h(Q) - \hat{h}(Q) \leq \frac{1}{4}h(j) + \frac{1}{6}h(\Delta) + 1.946.
\end{equation}
where $\Delta = -16(4a^3 + 27b^2)$ and $j = -48a^3/\Delta.$
\end{theorem}

We return to the proof of the final assertion in Theorem~\ref{upc}.
From~\eqref{defE}, we have
$$
C_n^2 = A_n^3 - 432m^2B_n^6.
$$
Suppose first that $|A_n|>B_n^2$. Since $\gcd(A_n,C_n) \mid 432m^2$, for any $\e>0$, the $ABC$ conjecture gives
$$
|A_n^3| \ll \left|mA_nB_nC_n\right|^{1+\e}
$$
and, since $|C_n|<|A_n|^{3/2}$, we get
$$
|A_n|^{1/2-\e} \ll mB_n.
$$
Thus, taking $\e=1/4$ and writing $|A_n|=h(nP)$, we have
\begin{equation}\label{BnClogm}
\log B_n >\tfrac 14 h(nP) - K_1\log m,
\end{equation}
for some constant $K_1>0$. On the other hand, if $|A_n|<B_n^2$ then~\eqref{BnClogm} is trivially satisfied, since $h(nP)=B_n^2$.

Finally, from Theorem~\ref{WeilCan} we have
$$
h(nP) > \hat h(nP) - K_2\log m = n^2 \hat h(P) - K_2\log m,
$$
for some constant $K_2>0$. Moreover, $\hat h(P) > K_3\log m$, for some constant~$K_3>0$, by~\cite[Proposition~1]{jed}.
Putting these together with~\eqref{BnClogm}, we get
$$
\log B_n > (K_4 n^2 - K_5)\log m,
$$
for constants $K_4,K_5>0$. In particular, this bounds, independently of $m$ and $P$, the index $n$ for which we may have $B_n=1$.
}
\end{proof}

\section{Proof of Theorem \ref{expupc}}\label{proof1.2}

The idea of the proof of Theorem~\ref{expupc} is to use heights on the Weierstrass elliptic
curve~\eqref{defE} to show that, for $n>22$, we are guaranteed to get a prime divisor of~$W_n$
coming from the term~$B_n$; that is, we use the isogeny to make explicit the constants in the
argument above using the $ABC$ conjecture. As we have seen, we also get a prime divisor coming
from~$A_n$ in this case, so that $W_n$ is not a prime power. The remaining values $2\le n\le 22$ are
treated case-by-case, using the explicit form of the point $nR$ as a rational function in
the coordinates $(u,v)$ of $R\in C(\mathbb Q)$.

\begin{lemma}\label{>m}
Let $E$ be defined as in \eqref{defE} and suppose $P\in E(\mathbb Q)$ is the image of a point in $E'(\mathbb Q)$
under the isogeny $\sigma$ as in \eqref{isog}. Then, with~$B_n$ as in \eqref{VnWn},
$$
B_n > 1
$$
for all $n > 22$.
\end{lemma}

\begin{proof}
Let $P \in E(\mathbb{Q})$ such that $\sigma(P')=P$, for some $P' \in E'(\mathbb Q)$. Write
$$
x_n := x(nP') = \dfrac{a_n}{b_n^2},
$$
with $\gcd(a_n, b_n) = 1$. From \eqref{isog},
\begin{equation}\label{Xn}
\dfrac{A_n}{B_{n}^{2}} = X(nP) = x_n+\frac{64m^2}{x_n^2} = \dfrac{a_{n}^{3}+64m^{2}b_{n}^{6}}{a_{n}^{2}b_{n}^{2}}.
\end{equation}

{\bf Claim:} $B_n > 1$,
provided $\max \{|a_n|,b_n^2\} > 8m$.

Before this claim can be settled, we must examine the fraction on the right-hand side of \eqref{Xn}
for possible cancellation. Let $d = p^r$ be a common factor of
$(a_{n}^{3}+64m^{2}b_{n}^{6})$ and $a_{n}^{2}b_{n}^{2}$, where $p$ is a prime and $r \in \mathbb{N}$
is the highest order of $p$ dividing both terms.
Since $\gcd(a_n, b_n) = 1$, either $d \mid  a_n^2$ or $d \mid  b_n^2$.
If the latter occurs, then $d\mid (a_{n}^{3}+64m^{2}b_{n}^{6})$ implies $d \mid a_n^3$, which is
impossible as $a_n$ and $b_n$ are coprime. Thus $d$ can only come from the term $a_n^2$, so
that $d \mid  a_n^3$. We have now that
$$
d\mid  (a_{n}^{3}+64m^{2}b_{n}^{6}),~d \mid a_n^3,\mbox{ and }d \nmid b_n^6,
$$
so $p^r = d \mid  64m^2.$ Hence the greatest common divisor of numerator and denominator of the
fraction on the right-hand side of \eqref{Xn}, say $g$, has to divide $64m^2$ as well.

To turn to the claim, if $|a_n| > 8m$, then
$$
B_n^{2} = \dfrac{a_n^{2}b_n^{2}}{g} \geq \dfrac{a_n^{2}}{g} \geq \dfrac{a_n^{2}}{64m^2} > \dfrac{64m^2}{64m^2}=1.
$$
On the other hand, if $b_n^2>8m$, then $g \mid a_n^2$ implies
$$
B_n^2 = \dfrac{a_n^{2}b_n^{2}}{g} \ge b_n^2 > 8m > 1.
$$
Thus we need to ensure that $\max \{ |a_n|, b_n^2\} > 8m,$
for our purposes and to that end we turn.
Note that the logarithm of the expression on the left is the Weil height $h(nP')$ of~$nP'$ which we know, from Theorem~\ref{WeilCan}, is close to the canonical height $\hat h(nP')$.

Write $h=\hat h(P)$ and $h'=\hat{h}(P')$; then
$$
h = \hat{h}(P) = \hat h(\sigma (P')) = 3 \hat{h}(P') = 3h'
$$
as $\sigma$ is a 3-isogeny. Applying the inequality~\eqref{diff} to the curve $E'$ with $Q=nP'$, $j=0$ and $\Delta = -16^{3}3^{3}m^4$, we obtain
\begin{equation}\label{bound}
\log \max\left\{|a_n|,b_{n}^2\right\} = h(nP') > h'n^2 -\frac{2}{3}\log{m} - \frac{1}{2}\log{48} - 2.14.
\end{equation}
The height bound \cite[(14)]{spag} gives
\begin{equation}\label{heightbound}
h' = \dfrac{h}{3} > \dfrac{1}{81}\log m - 0.039
\end{equation}
for all $m \geq 1$. Then~\eqref{bound} becomes
$$
\log \max\left\{|a_n|, b_{n}^2\right\} > \left(\dfrac{1}{81}\log m -0.039\right)n^2 - \dfrac{2}{3}\log{m} - \dfrac{1}{2}\log{48} - 2.14.
$$
Then we can ensure $\max \{|a_n|, b_n^2 \} > 8m$
provided that
\begin{equation}\label{overall}
\left(\dfrac{1}{81}\log m -0.039\right)n^2 - \dfrac{2}{3}\log{m} - \dfrac{1}{2}\log{48} - 2.14
 >\log (8m).
\end{equation}
With a manipulation, \eqref{overall} will be guaranteed for $n > 12$ but only for all sufficiently large $m$.
However the amount of checking for smaller values of~$m$ is infeasible.
More realistically, if $m > 353$ then \eqref{overall} is true provided $n > 22$.
Thus it can be concluded that for all $m>353$, we have
$B_n > 1$ if $n > 22$. The proof
for values $m\le 353$ follows in an appendix, see section~\ref{dontforget}.
\end{proof}

We are now in a position to prove Theorem \ref{expupc} using Lemma \ref{>m}.

\begin{proof}[Proof of Theorem \ref{expupc}]
When $n>12$, the term $A_n$ is not completely cancelled from the denominator of~\eqref{VnWn},
exactly as before. Also, it follows directly from Lemma~\ref{>m} that~$B_n>1$ for all~$n>22$.
Therefore, the term $W_n$ possesses at least two coprime factors for all $n>22$.
We will go on to prove the same for every $2\le n \leq 22$ case by case.

Before this, we give a proof of the second claim in Theorem~\ref{expupc}.
Under the isogeny hypothesis, it follows from \cite[Theorem 1.4]{eims} that
the number of prime power terms~$B_n$ is uniformly bounded. Combining this
with our knowledge that~$A_n$ has a primitive prime divisor for $n>12$ means that the
number of terms~$W_n$ which have at most two distinct prime divisors is uniformly
bounded.

With sleeves rolled up, we will now show that every term~$W_n$ fails
to be a prime power when $2<n\le 22$. We will not invoke the isogeny
hypothesis so these results apply in complete generality; that is, we are giving a proof of Lemma~\ref{lem:2<n<23}. Finally we
will deal with the case when $n=2$ assuming the isogeny hypothesis.

\begin{proof}[Proof of~Lemma~\ref{lem:2<n<23}]
Since $(W_n)$ is a divisibility sequence, it suffices to consider $W_n$ when $n$ is an odd prime less than $22$ or $n=4$. We begin with the case when $R=(u,v)$ is an integral point; at the end of the proof we will explain how the general case follows.

{\bf{Case $\bs n \bs= \bs4$}}
We deal firstly with the case when $n = 4$. Write
$$
4R = \left( \dfrac{U_4}{W_4}, \dfrac{V_4}{W_4} \right) =  \left( \dfrac{f_4(u,v)}{g_4(u,v)}, \dfrac{f'_4(u,v)}{g_4(u,v)} \right),
$$
where
$$
\dfrac{f_4(u,v)}{g_4(u,v)} = \dfrac{-u^{16} + 8v^3u^{13} + 32v^6u^{10} + 28v^9u^7 + 10v^{12}u^4 + 4v^{15}u }
{-u^{15} - 13v^3u^{12} - 10v^6u^9 + 10v^9u^6 + 13v^{12}u^3 + v^{15} },$$
and
$$
\dfrac{f'_4(u,v)}{g_4(u,v)} =  \dfrac{v^{16} - 8u^3v^{13} - 32u^6v^{10} - 28u^9v^7 - 10u^{12}v^4 - 4u^{15}v}
  {-u^{15} - 13v^3u^{12} - 10v^6u^9 + 10v^9u^6 + 13v^{12}u^3 + v^{15} }.
$$
We may consider the second coordinate, and factorize $g_4(u,v)$ as
a product of four terms:
\begin{eqnarray*}
g_{4,1}(u,v) &:=& v - u \\
g_{4,2}(u,v) &:=& u^2  + uv + v^2\ \equiv\ (v-u)^2 \pmod{3}\\
g_{4,3}(u,v) &:=& u^4  + 2u^3v + 2uv^3  + v^4\ \equiv\ (v-u)^4 \pmod{3}\\
g_{4,4}(u,v)  &:=& u^8  - 2u^7v + 4u^6v^2  + 4u^5v^3 - 5u^4v^4\\
&&\hskip2.5cm   + 4u^3v^5  + 4u^2v^6  - 2uv^7  + v^8 \\
&\equiv& (v-u)^8 \pmod{3}
\end{eqnarray*}
We claim that at least two of these factors can avoid being cancelled by the numerator $f'_4(u,v)$.
Choosing to consider $g_{4,3}$ and $g_{4,4}$, we can see that the
resultants between them and $f'_4$ with respect to $u$ and $v$ are
$$
R_u(f'_4, g_{4,3}) = 3^{16}v^{64}\mbox{ and }
R_v(f'_4, g_{4,3}) = 3^{16}u^{64},
$$
respectively, and also
$$
R_u(f'_4, g_{4,4}) = 3^{32}v^{128}\mbox{ and }
R_v(f'_4, g_{4,4}) = 3^{32}u^{128}.
$$
As $u$ and $v$ are coprime,
$$
\gcd(f'_4(u,v), g_{4,3}(u,v)) \mid  3^{16}\mbox{ and }\gcd(f'_4(u,v), g_{4,4}(u,v)) \mid  3^{32}.
$$
Next we will show that both $g_{4,3}(u,v)$ and $g_{4,4}(u,v)$ are not equal, up to a sign, to any power of $3$.
Suppose, for a contradiction, that $g_{4,3}(u,v) = \pm 3^k$,
for some $k > 1$. Then
$$
(v-u)^4  \equiv g_{4,3}(u,v) \equiv 0 \pmod{3}.
$$
Hence $u \equiv v \pmod{3}$, so $u^3 \equiv v^3 \pmod{3^2}$.
Replacing this in the expression of $g_{4,3}(u,v)$, we get
$$
0 \equiv u^4 + 2u^3v + 2u^4 + u^3v \equiv 3u^3(u+v) \pmod{3^2}.
$$
Then $3 \mid  u$ or $3 \mid  (u+v)$. Since $u \equiv v\pmod{3}$, the former implies $3 \mid  v$, and
the latter implies $3 \mid u $ and $3\mid v$ which is a contradiction as $\gcd(u,v) = 1$.
Thus the possibilities for $k$ such that $g_{4,3}(u,v) = \pm 3^k$ are only $0$ and $1$. Calculating
with PARI-GP~\cite{parigp}, without assuming the GRH,
shows that the only solutions $(u,v)$ of the equation $g_{4,3}(u,v) = \pm 1$ are $[0, \pm 1], [\pm 1, 0]$
and there are no solutions to $g_{4,3}(u,v) = \pm 3$.

A similar argument will be applied for the second factor $g_{4,4}(u,v).$ Suppose $g_{4,4}(u,v) = \pm 3^k$
for some $k > 2$. As $ (v-u)^8 \equiv g_{4,4}(u,v) \equiv 0  \pmod{3},$
we have $u \equiv v \pmod{3}$, so that
$$
u^3 \equiv v^3 \pmod{3^3},\quad 10u^3 \equiv v^3 \pmod{3^3},\quad\mbox{or }19u^3 \equiv v^3 \pmod{3^3}.
$$
Replacing each of these in the expression of $g_{4,4}$, we find that there are no solutions to
$g_{4,4}(u,v) = \pm 3^k$ when $k > 2$. Thus it remains to solve the equations $g_{4,4}(u,v) = \pm 3^k$ when $0 \leq k \leq 2$.
The only solutions to
$$
g_{4,4}(u,v) = \pm 1 \mbox{ are } [\pm 1, 0], [0, \pm 1], [-1,1], [1, -1].
$$
However these correspond to values of $m$ ($\pm 1$ or $0$) which yield rank zero (or singular) curves. Finally, GP says there
are no solutions at all to $g_{4,4}(u,v)=\pm 3^k$ when $k=1,2$.

We will prove moreover that the multiple $g_{4,3}(u,v)g_{4,4}(u,v)$ cannot be a prime power.
As above, $g_{4,3}$ and $g_{4,4}$ are not powers of $3$,
so write
$$
g_{4,3}(u,v) = \pm 3^{m}p_{1}^{m_1}\cdots p_{r}^{m_r}\mbox{ and }g
_{4,4}(u,v) = \pm 3^{n}q_{1}^{n_1}\cdots q_{s}^{n_s},
$$
where the $p_i$ and $q_j$ are primes other than $3$. Considering the resultant
between $g_{4,3}$ and $g_{4,4}$, we get $\gcd(g_{4,3}(u,v), g_{4,4}(u,v)) \mid  3^{10}$.
Thus there is at least one prime $p_i$ which is not equal to any prime $q_j$.
This implies $W_4$ is not a prime power.

{\bf{Case $\bs n \bs= \bs3$}} The expression for $3R$ can be written as
$$
3R = \left( \dfrac{u^9 + 6u^6v^3 + 3u^3v^6 - v^9}{3uv(u^6 + u^3v^3 + v^6)},
\dfrac{-u^9 + 3u^6v^3 + 6u^3v^6 + v^9}{3uv(u^6 + u^3v^3 + v^6)}  \right).
$$
For convenience, let
$$
f_3(u,v) = -u^9 + 3u^6v^3 + 6u^3v^6 + v^9\mbox{ and }g_3(u,v) = u^6 + u^3v^3 + v^6.
$$
By the theory of resultants, we obtain
$$
\gcd(f_3(u,v), g_3(u,v)) \mid  3^9.
$$
Since at least one of $u,v$ is not a unit, to complete the proof in this case
we have to prove that the denominator $g_3(u,v)$ is not (up to sign) a power of $3$.
Suppose $g_3(u,v) = \pm 3^k$ for some $k > 1$.
Then $(u-v)^6 \equiv g_3(u,v) \equiv 0\pmod{3}$.
Thus $u^3 \equiv v^3\pmod{3^2}$, and hence
$$
0 \equiv u^6 + u^3v^3 + v^6 \equiv 3u^6 \pmod{3^2},
$$
so $3\mid u$. This implies $3\mid v$ which is impossible.
For the remaining values, when $k=1$ the solutions are $[1, 1]$ and $[-1, -1]$ and
when $k = 0$, the only solutions are given by
$$
(u,v) = [\pm 1, 0], [0, \pm 1], [1, -1], [-1, 1].
$$
As before, these all correspond to impossible values of $m$ ($0,\pm1,\pm2$).
Since $\gcd(u,v) = 1$ and $u$ and $v$ are coprime to both $f_3(u,v)$ and $g_3(u,v)$ but are not
both units (as $m\ne 0,\pm 2$), $W_3$ possesses at least two coprime divisors.

{\bf{Case $\bs n \bs\equiv \bs 1\bs{\pmod{\bs 3}}$}} The proof in this case proceeds exactly
in the same way as in the case $n = 4$,
by the following steps.
\begin{list}{}{\leftmargin=0.8cm\itemsep=0.3cm\itemindent=-0.5cm}
\item[(i)] Write
$$
nR = \left( \dfrac{U_n}{W_n}, \dfrac{V_n}{W_n} \right) =  \left( \dfrac{f_n(u,v)}{g_n(u,v)}, \dfrac{f'_n(u,v)}{g_n(u,v)} \right),
$$
and factor the denominator $g_n(u,v)$ as $g_{n,1}(u,v), g_{n,2}(u,v),...,g_{n,k}(u,v)$, all of
which are homogeneous in $u$ and $v$.
By the theory of resultants, we have found fortunately that for each $n$,
$\gcd(f'_n(u,v), g_{n,i}(u,v))$ divides a power of $3$ for every $i = 1,...,k$.

\item[(ii)] Pick two factors of $g_n$, say $g_{n,i}(u,v)$ and $g_{n,j}(u,v)$.
It can be proved that both of them cannot
be a power of $3$ (up to sign) by using the following facts:
\begin{eqnarray*}
g_{n,i}(u,v) &\equiv& (u-v)^{\deg(g_{n,i})}\pmod 3,\\
g_{n,j}(u,v) &\equiv& (u-v)^{\deg(g_{n,j})}\pmod 3.
\end{eqnarray*}

\item[(iii)] Show that the multiple $g_{n,i}g_{n,j}$ is not a prime power, for
which it is sufficient to prove that the resultant of $g_{n,i}$ and $g_{n,j}$ is a power of~$3$.
\end{list}

{\bf{Case $\bs n \bs\equiv \bs 2\bs{\pmod{\bs 3}}$}} In this case, the situation is more complicated.
For all $n$, $f'_n(u,v)$ and $g_n(u,v)$
also behave as in the previous case in steps (i) and (iii). However, matters are slightly
different in step (ii).
We need to employ some facts about Newton polygons over $3$-adic fields to know about the
$3$-adic valuation of $g_{n,i}$.
We will show how to do this for $n = 5$ (for other $n$, the proofs will proceed in the same way).
We have
$$
g_{5,1}(u,v) = u^8-2u^7v-2u^6v^2+u^5v^3-5u^4v^4+u^3v^5-2u^2v^6-2uv^7+v^8,
$$
and
\begin{eqnarray*}
g_{5,2}(u,v) &\!\!\!=\!\!\!&
 u^{16}+2u^{15}v+6u^{14}v^2-2u^{13}v^3+11u^{12}v^4\\[1pt]
 &&\quad+21u^{11}v^5 -11u^{10}v^6-u^9v^7+27u^8v^8-
u^7v^9-11u^6v^{10}\\[1pt]
 &&\quad\ \quad\ +21u^5v^{11}
+11u^4v^{12}-2u^3v^{13}+6u^2v^{14}+2uv^{15}+v^{16}.
\end{eqnarray*}
We will explore their roots of the polynomials $h_{5,i}(X)=g_{5,i}(1+X,1)$, where we find
$$
h_{5,1}(X) = X^8 + 6X^7 + 12X^6 + 3X^5 - 30X^4 - 63X^3 - 63X^2 - 36X - 9,
$$
and
\begin{eqnarray*}
h_{5,2}(X) &\!\!\!=\!\!\!& X^{16} + 18X^{15} + 156X^{14} + 852X^{13} + 3261X^{12} + 9279X^{11}\\[1pt]
&&+ 20394X^{10} + 35496X^9 + 49617X^8 + 55971X^7 + 50814X^6\\[1pt]
&&\ \,+ 36774X^5 + 20871X^4 + 9072X^3 + 2916X^2 + 648X + 81.
\end{eqnarray*}
The Newton polygons for $h_{5,1}$ and $h_{5,2}$ with $p = 3$, as shown in Figures~\ref{g_51} and~\ref{g_52},
reveal that the slope of the only segment of each polygon is~$-\frac{1}{4}$.

\begin{figure}[!ht]
\begin{center}
\includegraphics[scale=0.65]{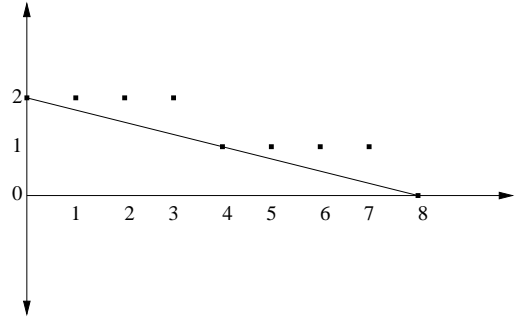}
\caption{Newton polygon of $h_{5,1}(X)$}
\label{g_51}
\end{center}
\end{figure}

By the standard facts about Newton Polygons, all roots of $h_{5,i}(X)$ (and all roots of $h_{11,i}(X)$ and $h_{17,i}(X)$) have
$3$-adic absolute values $3^{-\frac{1}{4}}$. If $\alpha$ is a root of $h_{5,i}(X)$ then, since $\left|u/v-1\right|_3 \ne |\alpha|_3$,
$$
\left| \frac uv - (1+\alpha) \right|_3
\geq |\alpha|_3 = 3^{-\frac{1}{4}}.
$$
Thus, since $|u|_3=|v|_3=1$,
$$
| g_{5,i}(u,v) |_3 = \left| h_{5,i}\left(\frac uv-1\right)\right|_3 =
\displaystyle{{\prod_{\alpha}}\left| \frac uv - (1+\alpha) \right|_3}  \geq (3^{-\frac{1}{4}})^{\deg(g_{5,i})},
$$
where $\alpha$ ranges over all roots of $h_{5,i}(X)$;
that is, the $3$-adic valuation of $g_{5,i}(u,v)$ is at most
$\deg(g_{5,i})/4$. It remains only to solve the
equations $g_{5,i}(u,v) = \pm 3^k$, with
$0 \leq k \leq \deg(g_{5,i})/4$. We find that the only solutions
to $g_{5,1}(u,v) = \pm 3^0$ are $[0, \pm 1], [\pm 1, 0]$ and there are no
solutions to $g_{5,1}(u,v) = \pm 3^k$ for other $k$. Similarly, the only solutions for $g_{5,2}$
arise for $g_{5,2}(u,v) = 3^0$, with $[\pm 1, 0], [0, \pm 1], [-1, 1], [1, -1]$,
and for $g_{5,2}(u,v) = 3^4$, with $[1, 1]$ and $[-1, -1]$. As before, these correspond to
inadmissible values of~$m$.

\begin{figure}[!ht]
\begin{center}
\includegraphics[scale=0.65]{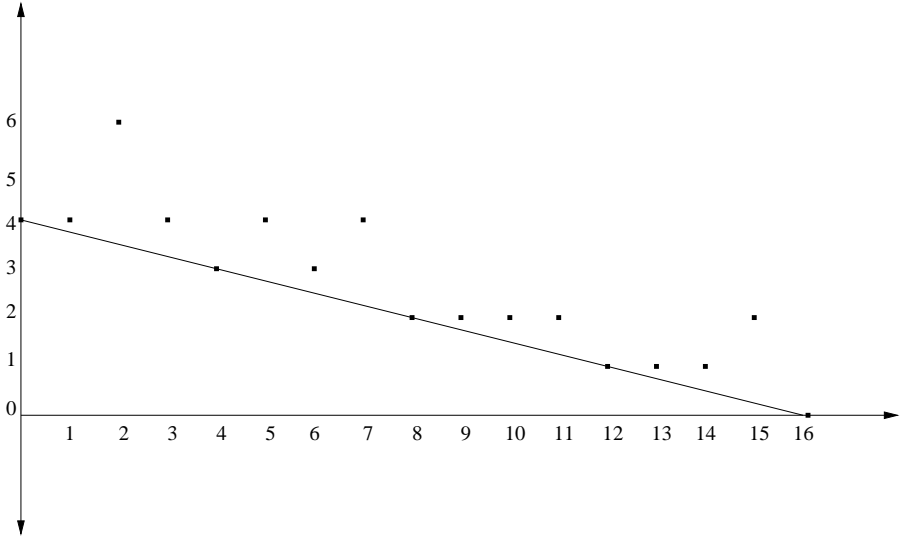}
\caption{Newton polygon of $h_{5,2}(X)$}
\label{g_52}
\end{center}
\end{figure}

This concludes the proof of Lemma~\ref{lem:2<n<23} when we consider integral points.
In the case of rational points, write
$$
R = \left( \dfrac{u_0}{w_0}, \dfrac{v_0}{w_0} \right) \in C(\mathbb{Q})
$$
in lowest terms.
The condition that $m$ is cube-free implies $u_0$ and $v_0$ are coprime.
Replacing $u$ and $v$ in the expressions of $nR$ in previous cases
by $u_0/w_0$ and $v_0/w_0$ respectively, we obtain
$$
nR = \left( \dfrac{f_n(u_0, v_0)}{w_0 g_n(u_0, v_0)}, \dfrac{f'_n(u_0, v_0)}{w_0 g_n(u_0, v_0)} \right).
$$
Now proceed with the proof working with $f_n(u_0, v_0)$ and $g_n(u_0, v_0)$: the conclusion follows as before.
Thus we have proved Lemma~\ref{lem:2<n<23}.
\end{proof}

{\bf{Case $\bs n \bs= \bs 2$}}
To complete the proof of Theorem~\ref{expupc} it remains only to treat the case $n=2$; the salient details follow.
Here we begin with
the curve~$C'$ in $(\ref{defC'})$. In a change of notation, we write $(u,v)\in C'(\mathbb Q)$
for a point which maps to~$R$ under the isogeny~$\tau$. For now we assume it is integral, noting that this
implies~$\gcd(u,v)=1$ since~$m$ is cube-free. We write
$$
2R = \left( \dfrac{f_2(u,v)}{g_2(u,v)}, \dfrac{f'_2(u,v)}{g_2(u,v)} \right),
$$
where
\begin{eqnarray*}
g_2(u,v)&\!\!=\!\!&-3(u-v)(u+2v)(2u+v)(u^2+uv+v^2)\\
&&\quad\ (u^6+3u^5v+60u^4v^2+115u^3v^3+60u^2v^4+3uv^5+v^6).
\end{eqnarray*}
Write $g_{2,1}(u,v)$ for the product of the linear factors and $g_{2,2}(u,v)$
for the degree~$6$ factor.

Firstly, we claim that the factor $g_{2,1}(u,v)$ does not completely cancel with the numerator $f_2(u,v)$. This is
easy to check: compute the resultant of each of the linear factors in $g_{2,1}(u,v)$ with $f_2(u,v)$.
Every time you obtain a power of~$3$. If $g_{2,1}(u,v)$
cancels then each of the linear factors is a power of~$3$ (possibly~$3^0$). This cannot happen since
$$
(2u+v)-(u+2v)=u-v
$$
and the equation $3^a+3^b=3^c$ has no solutions.

Similarly, the factor $g_{2,2}(u,v)$ is not cancelled by $f_2(u,v)$.
Checking resultants
shows any common division is a power of~$3$. Now $g_{2,2}(u,v)\equiv 0\pmod{3^6}$
forces $u\equiv v\equiv 0\pmod{3}$ so we only need to consider
$$
g_{2,2}(u,v)=\pm 3^r \mbox{ with } \gcd(u,v)=1 \mbox{ and } r=0,\dots ,5.
$$
The only solutions occur when $g_{2,2}(u,v)= 1$ and $g_{2,2}(u,v)=3^5$. In the first
case, they are
$$
[1, 0], [-1, 0], [1, -1], [0, -1], [-1, 1], [0, 1],
$$
each giving the inadmissible value $m=0$. In the second case, they are
$$
[1, -2], [2, -1], [-1, -1], [-2, 1], [1, 1], [-1, 2],
$$
giving the inadmissible values $m=\pm 2$.

Finally, the only common divisor of $g_{1,2}(u,v)$ and $g_{2,2}(u,v)$ is a power of~$3$
so we must obtain one non-trivial factor of $W_2$ from the first term and a coprime
factor from the second, which proves that~$W_2$ is not a prime power.

When $(u,v)$ is \emph{not} integral, things are slightly different from the situation for $2<n\le 22$ above,
because $u$ and $v$ might have different
denominators. Write $u=a/hb$ and $v=c/hd$ with $\gcd(a,b)=\gcd(c,d)=\gcd(b,d)=1$. Then also
$\gcd(a,c)=1$ (as $m$ is cube-free). Writing $u_0=ad$ and $v_0=bc$, which are coprime, we obtain
$$
2R = \left( \dfrac{f_2(u_0,v_0)}{bdh g_2(u_0,v_0)}, \dfrac{f'_2(u_0,v_0)}{bdh g_2(u_0,v_0)} \right)
$$
so we can proceed as before.
\end{proof}

\section{Appendix}\label{dontforget}

In order to complete the proof of Lemma \ref{>m}, it remains to check the
statement $(\ref{overall})$ for all cube-free integers $m$ up to $353$, as mentioned
at the end of the proof. In this part, we deal with the particular
computations to find a uniform bound, $N_0$, on the indices $n$ for such $m$.
Ranks and generators of $E: Y^2 = X^3 - 432m^2$
were computed using MAGMA \cite{magma}. Note that when $m = 337$, we
were unable to find the generator and rank
using MAGMA. This was found instead using SAGE \cite{sage}.
For rank-1 curves, we tested the elliptic divisibility sequence $(B_n)$
arising from the generator for $n = 1,...,22$.
A special argument is
required for the curves of rank 2, with two parts needed to find
the bound $N_0$.
Firstly find the finite set of pairs $(i,j)$, $i,j
\in \mathbb Z$, such that the canonical height of each point $iP + jQ$
is less than $40$,
where $P$ and $Q$ represent the generators. Then compute the
elliptic divisibility sequence $(B_n)$ arising from each point $iP + jQ$,
for $n = 1,...,22$. Now we get a bound, say $N_0'$, for the indices $n$
from the points of canonical height less than $40$. To treat all
cases, when $h > 40$, we return to the proof of Lemma \ref{>m} again, replacing the estimate
(\ref{heightbound}) by
\begin{center}
$h' > \dfrac{h}{3} > \dfrac{40}{3}.$
\end{center}
This leads to
\begin{center}
$\dfrac{40}{3}n^2 - \dfrac{2}{3}\log{m} - \dfrac{1}{2}\log{48} - 2.14
 > \log(8m).$
\end{center}
Taking specific values for $m$ such that $E$ has rank-$2$ gives another
bound, say $N_0''$, for the indices $n$. Comparing $N_0'$ and $N_0''$,
let
\begin{center}
$N_0 = \min\{N_0', N_0''\}.$
\end{center}
In no case did $N_0$ exceed $1$ for curves of rank $1$ and $2$.
There are no curves of higher rank appearing in that range.

\end{document}